\documentstyle[amscd,twoside]{amsart}

\input amssym.def
\input amssym

\newtheorem{theo}{\sc Theorem}[section]
\newtheorem{lemma}[theo]{\sc Lemma}
\newtheorem{corollary}[theo]{\sc Corollary}
\newtheorem{proposition}[theo]{\sc Proposition}
\newtheorem{definition}[theo]{\sc   Definition}

\title{$k$-very ample line bundles on Del Pezzo surfaces}
\author{Sandra Di Rocco}
\address{Sandra Di Rocco, Department of Mathematics\\ University of Notre
dame\\
 Mail Distribution Center\\ Notre Dame, Indiana 46556-5683}
\email{ sdirocco@@artin.helios.nd.edu}
\date{(revised copy)}

\begin{document}
\maketitle

\begin{abstract}
A $k$-very ample line bundle L on a Del Pezzo Surface is numerically
characterized, improving the results
of Biancofiore-Ceresa in [7].
\end{abstract}
\section{Introduction.}

Let L be a line bundle on a smooth connected projective surface. L is said
to be $k$-very ample
for an integer $k\geq 0$, if given any 0-dimensional subscheme $(Z,{\cal
O}_Z)$ of S with
$h^0({\cal O}_Z)=k+1$, the restriction map $H^0(L)\longrightarrow
H^0(L_{Z})$ is surjective. Notice
that L is spanned by its global section if and only if it is $0$-very ample
and that it is
very-ample if and only if it is 1-very ample. Using the Reider-type
criterion (theorem 2.3) it is possible to give an exact numerical
characterization of $k$-very ample line bundles on Surfaces whose Picard
group is fully
known. $k$-very ample line bundles on ${\bf P}^2$, on the Hirzebruch
surfaces ${\bf F_n}$ and on hyperelliptic surfaces
 are indeed completely characterized in [2], [5],
 [10].

 This paper gives a numerical
 characterization of $k$-very ample line bundles on Del Pezzo surfaces,
which improves the result
of Biancofiore and Ceresa in [7].

Let $S_r$ be a Del Pezzo surface of degree $9-r$, i.e. $-K_{S_r}$ is an
ample line bundle of degree
$9-r.$ We prove that:
\vspace{.20in}

{\em Let $L=al-\sum_1^rb_ie_i$ be a line bundle on $S_r$, $L\neq
-kK_{S_8}$, $L\neq -(k+1)K_{S_8}$ and
$L\neq -K_{S_7}$ when $k=1$, then $L$ is k-very-ample if and only if:\\
 for $r=1$ $a\geq b_1+k$ and $b_1\geq k$ ;\\
for $r=2,3,4$,  $b_1\geq b_2...\geq b_r\geq k$ and $a\geq b_i+b_j+k$, where
$i\neq j=1...r$;\\
for $r=5,6$, $b_1\geq ... \geq b_r\geq k$, $ a\geq b_i+b_j+k$, where $i\neq
j=1...r$, and
$ 2a\geq\sum_1^5b_{i_t}+k$;\\
for $r=7$, $b_1\geq ... \geq b_7\geq k$, $ a\geq b_i+b_j+k$, where $i\neq
j=1...r$,
$ 2a\geq\sum_1^5b_{i_t}+k$ and  $ 3a\geq 2b_i+\sum_1^6 b_{j_t}+ k$;\\
for $r=8$,  $b_1\geq ... \geq b_8\geq k$, $ a\geq b_i+b_j+k$ , $
2a\geq\sum_1^5b_{i_j}+k$ , \\
 $ 3a\geq 2b_i+\sum_1^6 b_{j_t}+ k$, $4a\geq\sum_1^32b_i+\sum_1^5b_{j_t}+k$,
 $5a\geq\sum_1^62b_{i_t}+b_j+b_k+k$
and  $6a\geq 3b_i+\sum_1^72b_{j_t}+k$}.
\newpage

We find that the set of the exceptional curves and effective divisors with
selfintersection
 zero, $\cal S$, plays a very important rule in checking the nefness and
$k$-very ampleness
of a line bundle on $S_r$.
 We show that a line bundle L is nef if and only if it is spanned and that L
is $k$-very ample if and only if the intersection with all the elements of
$\cal S$ is greater or
 equal than k. In the first part of the paper we give a complete numerical
description of the elements
of $\cal S$ that will allow us to prove the main theorem.\\

In [7] for L to be $k$-very ample it is required that $a\geq \sum_1^rb_i$
for any r. If we consider, for
example, the case $r=3$ Corollary 4.6 proves that it is enough to require
that $a\geq b_i+b_j$.\\

\noindent{\bf Acknowledgments.} I would like to thank my advisor, Andrew J.
Sommese, for suggesting this problem and for his help
throughout my studies. I would also like to thank Gianmario Besana for his
numerous helpful comments.
The author was supported by the C.N.R. of the Italian government.
\vspace{0.15in}
%%%%%%%%%%%%%%%%%%%%%%%%%%%%%%%%%%%%%%%%%%%%%%%%%%%%%%%%%%%%%%%
\section{Background material.}
%%%%%%%%%%%%%%%%%%%%%%%%%%%%%%%%%%%%%%%%%%%%%%%%%%%%%%%%%%%%%%%
\vspace{.18in}
\noindent{ \bf \sf  Notation.}
\vspace{.15in}

Let S be a smooth connected projected surface defined over the complex
field {\bf C}.
$B_{P_1,...,P_r}(S)$ is the blowup of S at the points $P_1,...,P_r$;\\
Let L be a line bundle on $S$, we will use the following notation:\\
$h^i(L)$ is the dimension of $H^i(L)$;\\
$|L|$ is the complete linear system associated to L;\\
$K_S$ is the canonical bundle of S;\\
$L^2=d(L)$ is the degree of the line bundle L;\\
$g(L)$, the sectional genus of S, is the integer defined by the equality :\\
$2g(L)-2=L^2+K_SL$. If $L={\cal O}_S(D)$ for a divisor D we also write
$g(D)$ for $g(L)$ and if D is an
irreducible reduced curve on a surface S , $g(D)$ is the arithmetic genus of D.
\\ \\
L is said to be {\it numerically effective} (nef) if $LC\geq 0$ for every
curve $C$ on S and in this case
L is said to be {\it big} if $L^2>0$.\\
L is said to be Q-effective if, for some positive integer $n$, $nD$ is
effective.\\
L is said to be spanned if it is spanned by its global sections, i.e.\\
$h^0(L\bigotimes{\cal I}_P)=h^0(L)-1$ for every point on S.\\
L is said to be very ample if the map associated to $|L|$ embeds S in
$P^{h^0(L)-1}$, i.e.
 $h^0(L\bigotimes{\cal I}_Z)
=h^0(L)-2$ for any $0$-dimensional subscheme Z of length 2 on S.\\
\newpage
A $(-1)$-curve $C$ on $S$, i.e. a rational curve of degree $-1$, will be
always called an exceptional
 curve.
\vspace{.15in}
%%%%%%%%%%%%%%%%%%%%%%%%%%%%%%%%%%%%%%%%%%%%%%%%%%%%%%%%%%%%%%%

\noindent{ \bf \sf  $k$-very ampleness. }
%%%%%%%%%%%%%%%%%%%%%%%%%%%%%%%%%%%%%%%%%%%%%%%%%%%%%%%%%%%%%%%
\begin{definition}
A line bundle is said to be {\bf k-very ample},
for an integer $k\geq 0$ if, given any 0-dimensional subscheme $(Z,{\cal
O}_Z)$ of S with
$h^0({\cal O}_Z)=k+1$, the restriction map $H^0(L)\longrightarrow
H^0(L_{Z})$ is surjective.
\end{definition}
Note that L is 0-very ample if and only if it is spanned by its global
section and L is 1-very ample if and only if
it is very ample.\\
The notion of $k$-very ampleness has a classical interpretation in terms of
associated secant mappings.
Let $S^{[k]}$ be the Hilbert scheme of 0-dimensional subschemes on S of
length k; a line bundle
L is $k$-very ample if an associated line bundle $L^{[k]}$, on $S^{[k]}$,
is very ample, see [3].\\
\begin{proposition}(see [4], 0.5.1) If L is a $k$-very ample line bundle on
a surface S, then $LC\geq k$ for
every effective curve C on S, with equality only if $C\cong {\bf P}^1$.
\end{proposition}
The following theorem will be the tool we will use to check the $k$-very
ampleness for line bundles:\\
\begin{theo}(see [3], Theorem 2.1)
Let M be a nef line bundle on a smooth surface S, such that $M^2\geq 4k+5$,
 for $k\geq 0$. Then either $L=M+K_S$
is $k$-very ample or there exists an effective divisor D s.t M-2D is
Q-effective, D contains some
$0$-dimensional subscheme where the $k$-very ampleness fails and:
$$ MD-k-1\leq D^2<MD/2<k+1$$
\end{theo}
\vspace{.25in}
%%%%%%%%%%%%%%%%%%%%%%%%%%%%%%%%%%%%%%%%%%%%%%%%%%%%%%%%%%%%%%%
\noindent{ \bf \sf  Del Pezzo surfaces.}
%%%%%%%%%%%%%%%%%%%%%%%%%%%%%%%%%%%%%%%%%%%%%%%%%%%%%%%%%%%%%%%
\vspace{.15in}

A surface is called a Del Pezzo surface if its anticanonical bundle $-K_S$
is ample.
If $K_S^2=s$ we will say that S is a Del Pezzo surface of degree s.
The degree turns out to be a very useful tool to classify Del  Pezzo
surfaces. Indeed we have
the following classification:
\begin{proposition}(see [8])
Let S, s as before, then\\
a) $1\leq s\leq 9$;\\
b) $s=9$ if and only if $S=P^2$;\\
c) if s=8 then either $S=P^1\times P^1$ or $S=B_P(P^2)={\bf F_1}$;\\
d) if $1\leq r\leq 7$ then $S=B_{P_1,...,P_{9-s}}(P^2)$ where no three points
 are on a line nor six points on a conic.
\end{proposition}
We will denote a Del Pezzo surface of degree s with $S_r$, where $r=9-s$.\\
Let $\pi :S_r\longrightarrow {\bf P}^2$ be the blowing up , as in d). It
follows that $Pic(S_r)={\bf Z}^{r+1}$ and it is generated by $\{l,e_1,...,
e_r\}$ where $l$ is the class of $\pi^*({\cal O}_{\bf  P^2}(1))$, $e_i$ the
class of $\pi^{-1}(P_i)$, $l^2=1$, $e_1^2=-1$, $e_1e_j=0$ and $e_il=0$.\\
Therefore any line bundle on $S_r$ is of the form $al-\sum_1^rb_ie_i$ and
from the definition of blow up
it follows that $-K_S=3l-\sum_1^re_i$.\\
Since for $S_0$ and for $S_1=P^1\times P^1$ a complete characterization of
$k$-very ample line bundles
has been done in [2] and [5], from now on we'll focus our attention on
$S_r=B_{P_1,...,P_r}$ for
$r=1,...,8$, as in proposition 2.4 d).\\
We will use the following notation, according to [8]:\\
Let $I_r=\{L\in Pic(S_r)$ such that $L^2=-1$ and $LK_S=-1\}$ be the set of
exceptional curves
 on $S_r$,
they are in fact irreducible effective divisors.
the cardinality of $I_r$ is finite for $1\leq r \leq 8$ and any $\xi\in
I_r$ can be expressed
as the pull back of a curve in {\bf P}$^2$ passing through the points blown up.
 In the following table the number of exceptional curves in $S_r$ of all
possible
 types it is shown
. For example if $r=6$, i.e $S_r$ is the cubic surface in $P^3$, the
table shows the 27 possible lines.
\vspace{.05in}
\begin{center}
\begin{tabular}{|l|} \hline
$type\backslash r$\\ \hline \hline
$(0;-1)$\\  \hline
$(1;1^2)$\\  \hline
$(2;1^5)$\\  \hline
$(3;2,1^6)$\\ \hline
$(4;2^3,1^5)$\\ \hline
$(5;2^6,1^2)$ \\ \hline
$(6;3,2^7)$ \\ \hline
\end{tabular}\hspace{.005pt}
\begin{tabular}{|l|l|l|l|l|l|l|l|l|} \hline
 1 & 2 & 3 & 4 & 5 & 6 & 7 & 8 \\ \hline\hline
 1 & 2 & 3 & 4 & 5 & 6 & 7 & 8  \\ \hline
 0 & 1 & 3 & 6 & 10 & 15 & 21 & 28 \\ \hline
 0 & 0 & 0 & 0 & 1 & 6 & 21 & 56 \\ \hline
 0 & 0 & 0 & 0 & 0 & 0 & 7 & 56 \\ \hline
  0 & 0 & 0 & 0 & 0 & 0 & 0 & 56 \\ \hline
 0 & 0 & 0 & 0 & 0 & 0 & 0 & 28 \\ \hline
 0 & 0 & 0 & 0 & 0 & 0 & 0 & 8 \\ \hline
\end{tabular}
\end{center}
\vspace{.04in}
where $(a_0;a_1^{n_1},a_2^{n_2},...)$ is the class of the curve
$a_0l-\sum_1^{n_1}a_1e_{i_t}-\sum_1^{n_2}
a_2e_{j_t}-...$, i.e the pull back of a curve of degree $a_o$ in ${\bf
P}^2$ passing,  with multiplicity
$a_i$ , through $n_i$  points blown up. \\
\begin {lemma}(see [8], Lemme 9)If D is an effective irreducible divisor on
$S_r$ such that
 $DK_S=-1$ then D is an exceptional curve or $D=-K_{S_8}$.
\end{lemma}
%%%%%%%%%%%%%%%%%%%%%%%%%%%%%%%%%%%%%%%%%%%%%%%%%%%%%%%%%%%%%%%
\section{ Nefness on Del Pezzo surfaces}
%%%%%%%%%%%%%%%%%%%%%%%%%%%%%%%%%%%%%%%%%%%%%%%%%%%%%%%%%%%%%%%
Let $D=al-\sum_1^rb_ie_i\in Pic(S_r)$ s.t $D^2=0$, i.e.  $a^2-\sum_1^rb_i^2=0$.
By the adjunction formula we have $g(D)=0$ and $K_SD=-2$ i.e.,
$3a-\sum_1^rb_i=2$. Using the inequality
$(\sum_1^rb_i)^2\leq r\sum_1^rb_i^2$ and:
$$\left\{\begin{array}{c}
\sum_1^r b_i  = 3a-2 \\
\sum_1^r b_i^2  = a^2
\end{array}
\right. $$
we get $a\leq 11$.\\
 Assuming $b_1\leq b_2\leq ... \leq b_r$, by numerical computations we
obtain that there is
 a finite number of n+1-tuples $(a,b_1,b_2,...,b_r)$ , solutions of the
system above. In the
following table we collect the possible solutions and show the divisors
obtained:
\vspace{0.05in}
\begin{center}
\begin{tabular}{|l|l|l|l|l|l|l|l|l|l|} \hline
$a$ &$b_1$ & $b_2$ & $b_3$ & $b_4$ &$ b_5$ &$b_6$ & $b_7$ & $b_8$ &type  \\
\hline\hline
1& 0& 0& 0 & 0 & 0 & 0 & 0 & 1&(1;1) \\ \hline
2& 0 & 0 & 0 & 0 & 1 & 1 & 1 & 1&$(2;1^4)$ \\ \hline
3& 0 & 0 & 1 & 1 & 1 & 1 & 1 & 2&$(3;1^5;2)$ \\ \hline
4 & 1 & 1 & 1 & 1 & 1 & 1 & 1 & 3&$(4;1^7;3)$ \\
  & 0 & 1 & 1 & 1 & 1 & 2 & 2 & 2&$(4;1^4;2^3)$  \\ \hline
5 & 0 & 1 & 2 & 2 & 2 & 2 & 2 & 2&$(5;1;2^6)$ \\
  & 1 & 1 & 1 & 1 & 2 & 2 & 2 & 3&$(5;1^4;2^3;3)$ \\ \hline
6 & 1 & 1 & 2 & 2 & 2 & 2 & 3 & 3&$(6;1^2;2^4;3^2)$ \\ \hline
7 & 1 & 2 & 2 & 2 & 3 & 3 & 3 & 3&$(7;1;2^3;3^4)$ \\
  & 2 & 2 & 2 & 2 & 2 & 2 & 3 & 4&$(7;2^6;3;4)$ \\ \hline
8 & 1 & 3 & 3 & 3 & 3 & 3 & 3 & 3&$(8;1;3^7)$ \\
  & 2 & 2 & 2 & 3 & 3 & 3 & 3 & 4&$(8;2^3;3^4,4)$ \\ \hline
9 & 2 & 3 & 3 & 3 & 3 & 3 & 4 & 4&$(9;2,3^5;4^2)$ \\ \hline
10 & 3 & 3 & 3 & 3 & 4 & 4 & 4 & 4&$(10;3^4,4^4)$ \\ \hline
11 & 3 & 4 & 4 & 4 & 4 & 4 & 4 & 4&$(11;3;4^7)$ \\ \hline
\end {tabular}
\end{center}
\vspace{.04in}
where every row represents a solution.

Note that all the divisors obtained are sum of two exceptional divisors,
except for the case $ a=1 $.
The decomposition, for any value of $a$, it is shown in the next table :\\
\begin{center}
\begin{tabular}{|l|l|}\hline
a & type \\ \hline\hline
1 & $(1;1)$ \\ \hline
2 & $(1;1^2)$+$(1;1^2)$ \\ \hline
3 & $ (2;1^5)+(1;1^2)$ \\ \hline
4 & $ (3;2,1^6)+(1;1^2)$ \\ \hline
5 & $ (3;2,1^6)+(2;1^5)$ \\ \hline
6 & $(4;2^3,1^5)+(2;1^5)$ \\ \hline
7 & $(5;2^6,1^2)+(2;1^5)$ or $(6;3,2^7)+(1;1^2)$ \\ \hline
8 & $(5;2^6,1^2)+(3;2,1^6)$ or $(6;3,2^7)+(2;1^5)$ \\ \hline
9 & $(6;3,2^7)+(3;2,1^6)$ \\ \hline
10 & $(6;3,2^7)+(4;2^3,1^5)$ \\ \hline
11 & $(6;3,2^7)+(5;2^6,1^2)$ \\ \hline
\end{tabular}
\end{center}
\vspace{.08in}
The following proposition follows from the computations above:
\begin{proposition}Let $D$ be an effective divisor on $S_r$ with $D^2=0$,
then:  \\
a) D is irreducible and $D=l-e_i$; \\
b) D is reducible and $D=\xi_1+\xi_2$, $\xi_1,\xi_2\in I_r$.
\end{proposition}
\begin{proposition}Let D be an effective divisor on $S_r$, then it is nef
and big, unless it is irreducible
and $D\in I_r$ or $D=l-e_i$, or it is reducible and it has at least one
exceptional curve as component.
\end{proposition}
\begin{pf}  If S is a surface such that for every two divisors $A$, $B$,
with $A^2\geq0$ and $B^2
\geq 0$, $AB\geq 0$, then the same property is true for a surface obtained
by blowing up one
point. If $S$ is ${\bf P}^2$
the property is true, therefore it is true for $S_r$.
The only effective divisors on $S_r$ with negative self-intersection are
 the elements of $I_r$. Thus If $D$ is an irreducible effective divisor, $
D^2>0$ unless
$D\in I_r$ or $D=l-e_i$ and $D\xi \geq 0$ for any $\xi\in I_r$, unless
$D=\xi\in I_r$.
 If D is reducible then $D^2>0$ unless $D=\xi_1+\xi_2$ with $\xi_1,\xi_2\in
I_r$
 and  if $D\xi<0$ for some $\xi\in I_r$ then $\xi$ is one of its components.
\end{pf}
\begin{lemma}(see [8], Lemme 1) Let $L=al-\sum_i^rb_ie_i\in Pic(S_r)$. If
$a\geq -2$ and $d(L)\geq K_SL$
 then L is effective.
\end{lemma}
\begin{theo} Let $L=al-\sum_1^3b_ie_i$ be a line bundle on $S_r$.
 L is nef if and only if $L\xi\geq 0$ for any $\xi\in I_r$ and
$L(l-e_1)\geq 0$ for $r=1$ i.e.\\
 for $r=1$, $a\geq b_1$ and $b_1\geq 0$ ;\\
for $r=2,3,4$,  $b_1\geq b_2...\geq b_r\geq 0$ and $a\geq b_i+b_j$;\\
for $r=5,6$, $b_1\geq ... \geq b_r\geq 0$, $ a\geq b_i+b_j$, where $i\neq
j=1...r$ and
$ 2a\geq\sum_1^5b_{i_t}$;\\
for $r=7$, $b_1\geq ... \geq b_7\geq 0$, $ a\geq b_i+b_j$, where $i\neq
j=1...r$,
$ 2a\geq\sum_1^5b_{i_t}$ and $ 3a\geq 2b_i+\sum_1^6 b_{j_t}$;\\
for $r=8$,  $b_1\geq ... \geq b_8\geq 0$, $ a\geq b_i+b_j$, where $i\neq
j=1...r$,
$ 2a\geq\sum_1^5b_{i_t}$, $ 3a\geq 2b_i+\sum_1^6 b_{j_t}$,
$4a\geq\sum_1^32b_{i_t}+\sum_1^5b_{j_t}$,
 $5a\geq\sum_1^62b_{i_t}+b_{j_t}+b_k$ and \\ $6a\geq
3b_{i_t}+\sum_1^72b_{j_t}$.
 \end{theo}
\begin{pf} If L is nef $L\xi\geq 0$ by definition.\\
Let L be a line bundle on $S_r$ s.t. $L\xi\geq 0$ for any $\xi\in I_r$ and
$L(l-e_1)\geq 0$ for $r=1$, then $a\geq 0$, $L^2\geq 0$ and $K_SL<L^2$ that
implies that L is effective on $S_r$.
Let $C$ be an effective divisor then, by proposition 3.2, $C$ is nef and
therefore $LC\geq 0$ unless
$C\in I_r$, $C=l-e_i$ or $C$ has an exceptional curve as component. If
$C\in I_r$
 or  has an exceptional curve as component $LC\geq 0$ by hypothesis.
If $C=l-e_i$ and $LC< 0$ then $r\geq 2$ and $a<b_i$; it follows that
$L\xi<0$ for $\xi=l-e_i-e_j$ that is impossible by
hypothesis.\end{pf}
%%%%%%%%%%%%%%%%%%%%%%%%%%%%%%%%%%%%%%%%%%%%%%%%%%%%%%%%%%%%%%%
\section{ $k$-very ampleness on Del Pezzo surfaces.}
%%%%%%%%%%%%%%%%%%%%%%%%%%%%%%%%%%%%%%%%%%%%%%%%%%%%%%%%%%%%%%%
Note that if $L=al-\sum_1^rb_ie_i$ is a $k$-very ample line bundle on $S_r$
then $LC\geq k$ for any
irreducible effective divisor C ; in particular $Ll\geq k$, $Le_i\geq k$
for $i=1...r$ and
$L\xi\geq k$ for any $\xi$ in $I_r$.\\
Let us assume $L\xi\geq k$ for any $\xi\in I_r$ and $l(l-e_i)\geq k$; that
implies L nef and big
 $b_i\geq k$ for $i=1...r$ and
$a\geq 3k$, for $r\geq 2$. Therefore we can write $L=kL'+H$ where $L'$ is
an effective divisor with
all coefficients greater or equal then k and $H$ is an effective divisor
with all coefficients strictly
less than k.
\begin{lemma} Let $M=L-K_S$; if $(M-D)D\leq k+1$ for some effective divisor
D on S, then \\
a) $r\geq 2$ and either $D^2>0$ or $L\xi <k$ for some $\xi\in I_r$;\\
b) $r=1$ and either $D^2>0$ or $Le_1< k$  or $L(l-e_1)< k$.
\end{lemma}
\begin{pf}  Assume $D^2<0$, then $D=\xi\in I_r$, $MD\leq k$ and $-K_SD\geq
1$, therefore $LD<k$.
Assume now $D^2=0$; if D is irreducible it is of the form $l-e_i$ with
$-K_SD=2$, that implies
$LD<k$. For $r\geq 2$, taking any exceptional curve of the form
$\xi=l-e_i-e_j$, we obtain $L\xi<k$.
If D is reducible then $r\geq 2$ and $D=\xi_1+\xi_2$ $\xi_1,\xi_2\in I_r$
which gives $L\xi_i<k$,
 for i=1,2.
\end{pf}
\begin{lemma} Let $L=kL'+H$  be a line bundle  on $S_r$ such that $L\xi\geq
k$ for any exceptional
 curve $\xi\in I_r$, and $L(l-e_1)\geq k$ for $r=1$, and $M=L-K_{S_r}$,
then either $(L')^2>0$ or $r=8$, $L=k\xi+H$ or $L=k(\xi_i+
\xi_j)+H$ with  $\xi,\xi_i,\xi_j$ of the
type $(5;2^6,1^2)$ or $(6;3,2^7)$ and there is no effective divisor $D$ on
$S_8$
such that $$ MD-k-1\leq D^2<MD/2<k+1$$.
\end{lemma}
\begin{pf} Assume $L'$ irreducible and $(L')^2\leq 0$. If $(L')^2=0$ then
$L'=l-e_i$ and therefore
$L(l-e_i-e_j)=H(l-e_i-e_j)<k$; since the coefficients of $H$ are strictly
less than k, that is
 impossible by our setting. If $(L')^2<0$ then $L=k\xi+H$ where $\xi\in
I_r$ and for $r\leq 7$ we
always have $\xi(l-e_1-e_2)=0$ for any $\xi\in I_r$ and therefore
$L(l-e_1-e_2)=H(l-e_1-e_2)<k$,
 impossible as before.
Assume now $r=8$ and $L=k\xi+H$  with  $\xi$ of the type $(5;2^6,1^2)$ or
$(6;3,2^7)$,
then $H^2>0$ and $H\xi\leq 2k$ otherwise $L\xi<k$.
If $D$ were an effective divisor on $S_8$ s.t. $ MD-k-1\leq D^2<MD/2<k+1$
then $D^2>0$, by
lemma 4.1.
 If $D^2=1$ then, since $(k\xi+H)^2\geq-k^2+1+4k^2=3k^2+1$ and $L^2D^2\leq
(LD)^2$, $LD>k$ and
$LD-K_SD-k-1>1=D^2$. If $D^2\geq 2$ then $LD>2k$, that implies
$k+1>D^2>2k+2-k-1$.
The case L irreducible, i.e $L=\xi_1+\xi_2$ with $\xi$ of the
type $(5;2^6,1^2)$ or $(6;3,2^7)$ is proved as before.\end{pf}

\begin{lemma}$M^2\geq 4k+5$ for any k unless $L=-kK_{S_8}$ or $k=1$ and
$L=-K_{S_7}$.
\end{lemma}\begin{pf}  by Lemma 4.2 $(L')^2\geq 1$ and $L'H\geq 0$; note that
we
 can also assume $-K_SL\geq 2$ unless r=8 and $L=-K_{S_8}$ in which case
 $M^2\geq 4k+5$ for any k, unless $H\equiv 0$ and therefore $L=-kK_{S_8}$.
A numerical checking also
shows that $(-(k+1)K_{S_7})^2\leq 4k+4$ if and only if k=1. In all the
remaining cases we have
$M^2\geq k^2+6k+2\geq 4k+5 $.\end{pf}

Note that if L is a k-very-ample line bundle on a surface then $LE\geq k+2$
for any elliptic curve
 E on S, therefore $-kK_{S_8}$ , $-(k+1)K_{S_8}$  cannot be $k$-very ample
and $-K_{S_7}$ cannot
be very ample.\\

\begin{theo} Let L be a line bundle  on $S_r$ such that $L\xi\geq k$ for
any exceptional
 curve $\xi\in I_r$, and $L(l-e_1)\geq k$ for $r=1$. L is
k-very-ample unless $L=-kK_{S_8}$, $L=-(k+1)K_{S_8}$ or $k=1$ and $L=-K_{S_7}$.
\end{theo}
\begin{pf} $M^2\geq 4k+5$ by Lemma 4.3  and  L is nef and big by theorem
3.4, then L is k-very
 ample unless there exist an effective divisor D on S such that $$
MD-k-1\leq D^2<MD/2<k+1.$$
 Assume such D exists.
 $(L')^2\geq 1$ by Lemma 4.2 and $D^2\geq 1$ by Lemma 4.1.
Therefore  $L'D\geq 1$ and $-K_SL'\geq 2$ unless r=8, $L'=-K_{S_8}$
and $H$ is not 0, in which case we would get $MD\geq 2k+2$ that is
impossible. If $L'D=1$ then by the Hodge index theorem we should have
$D^2=(L')^2=1$ and $-K_SD\leq 2$. If $-K_SD=1$ then $D=K_{S_8}$ that is
impossible since $-L'K_S\geq 2$.
 If $-K_SD=2$ we have $HD\leq 0$, i.e. $H^2\leq 0$ and $L^2\geq 3k^2-1$,
therefore $LD>2k$ and
$2k+2-k-1\leq 1$.
 If $L'D\geq 2$ then $MD\geq 2k+2$, that is impossible.\end{pf}
\begin{corollary}  Let L  be a line bundle on $S_r$, s.t $L\neq -kK_{S_8}$,\\
$L\neq -(k+1)K_{S_8}$ and
$L\neq -K_{S_7}$ when $k=1$. L is $k$-very ample if and only if $L\xi\geq
k$, for any $\xi\in I_r$
and $L(l-e_1)\geq k$ for r=1.
\end{corollary}
 Using the table of exceptional curves we easily obtain the following
numerical characterization:
\begin{corollary}  Let $L=al-\sum_1^rb_ie_i$ be a line bundle on
$S_r$,$L\neq -kK_{S_8}$, \\$L\neq -(k+1)K_{S_8}$ and
$L\neq -K_{S_7}$ when $k=1$ then $L$ is k-very-ample if and only if:\\
 for $r=1$, $a\geq b_1+k$ and $b_1\geq k$ ;\\
for $r=2,3,4$,  $b_1\geq b_2...\geq b_r\geq k$ and $a\geq b_i+b_j+k$, where
$i\neq j=1...r$;\\
for $r=5,6$, $b_1\geq ... \geq b_r\geq k$, $ a\geq b_i+b_j+k$, where $i\neq
j=1...r$, and
$ 2a\geq\sum_1^5b_{i_t}+k$;\\
for $r=7$, $b_1\geq ... \geq b_7\geq k$, $ a\geq b_i+b_j+k$, where $i\neq
j=1...r$,
$ 2a\geq\sum_1^5b_{i_t}+k$ and  $ 3a\geq 2b_i+\sum_1^6 b_{j_t}+ k$;\\
for $r=8$, $b_1\geq ... \geq b_8\geq k$, $ a\geq b_i+b_j+k$ ,$
2a\geq\sum_1^5b_{i_j}+k$  \\
 $ 3a\geq 2b_i+\sum_1^6 b_{j_t}+ k$, $4a\geq\sum_1^32b_i+\sum_1^5b_{j_t}+k$,
 $5a\geq\sum_1^62b_{i_t}+b_j+b_k+k$
and  $6a\geq 3b_i+\sum_1^72b_{j_t}+k$.
\end{corollary}

Comparing Theorem 3.4 and Corollary 4.6 we can conclude that:
\begin{theo}A line bundle on a Del Pezzo surface is nef if and only if it
is spanned by its global sections.
\end{theo}

\section{observations}

{\sf \bf The Hirzebruch surface ${\bf F_1}$.} For $S_1={\bf F}_1$ a
different characterization is
 often used, i.e. a line bundle L  is
 written as  $L=a_0E_0+bf$, where $E_0$ and $f$ are equal respectively to
the minimal section and
the fiber of the morphism $S_1\longrightarrow {\bf P}^1$ and $ K_S=-2E_0
-3f$. A simple checking
shows that $E_0=e_1$ and $f=l-e_1$.\\
 Therefore If $L\in Pic({\bf F}_1)$ is written as $a_0E_0+bf$
 we can express it as $L=bl-(b-a_0)e_1$.
{}From Corollary 2.5 it follows that L is $k$-very ample if and only if
$a_0\geq k$ and $b\geq a_0+k$, same characterization as in
[5].\\
\\
{\bf \sf The adjoint bundle.} Let $L$ be a $k$-very ample line bundle on
$S_r$. Applying
Corollary 4.6 to the adjoint bundle $K_S+L=(a-3)-\sum_1^r(b_i-1)e_i$ we can
easily see that it is always $(k-1)$-very ample if
$r\geq 2$.

If $r=1$ it is $(k-1)$-very ample if $a\geq b_1+k+1$.\\
\\
{\bf\sf The degree.} If L is $k$-very ample, from [5, theorem 3.1] $
kK_S+L$ is nef and big unless $L=-kK_S$, with $K_S^2\geq 2$
 for $k\geq 2$ and with $K_S^2\geq 3$ for $k=1$.

 From [1, Corollary 3.31] it follows that  $(kK_S+L)L\geq k+4$ for
$k\geq 2$ i.e. $$d(L)\geq k^2+3k+2$$ if $L\neq -kK_S$ and $k\geq 2$.

\bibliographystyle{abbrv}

\end{document}